\documentclass[12pt]{article}
\usepackage{amsthm,amsfonts,amssymb,amscd}

\begin{document}


\begin{center}
\large \bf The $4n^2$-inequality for\\
complete intersection singularities
\end{center}\vspace{0.3cm}

\centerline{A.V.Pukhlikov}\vspace{0.3cm}

\parshape=1
3cm 10cm \noindent {\small \quad\quad\quad \quad\quad\quad\quad
\quad\quad\quad {\bf }\newline The famous $4n^2$-inequality is
extended to generic complete intersection singularities: it is
shown that the multiplicity of the self-intersection of a mobile
linear system with a maximal singularity is greater than
$4n^2\mu$, where $\mu$ is the multiplicity of the singular point.

Bibliography: 16 titles.} \vspace{0.3cm}

\vspace{0.3cm}

\begin{flushright}
{\it To Askol'd Georgievich Khovanskii}
\end{flushright}

{\bf 1. Statement of the result.} Let $(X,o)$ be a germ of a
complete intersection singularity of codimension $l$ and type
$\underline{\mu}=(\mu_1,\dots,\mu_l)$, where
$$
\mathop{\rm dim} X=M\geqslant l+\mu_1+\dots +\mu_l+3.
$$
We will assume the singularity to be generic in the sense of Sec.
2 below. The aim of this note is to prove the following
claim.\vspace{0.1cm}

{\bf Theorem.} {\it Let $\Sigma$ be a mobile linear system on X.
Assume that for some positive $n\in{\mathbb Q}$ the pair
$(X,\frac{1}{n}\Sigma)$ is not canonical at the point $o$ but
canonical outside this point. Then the self-intersection
$Z=(D_1\circ D_2)$ of the system $\Sigma$ satisfies the
inequality}
\begin{equation}\label{24.06.2016.1}
\mathop{\rm mult}\nolimits_oZ>4n^2\mathop{\rm mult}\nolimits_oX.
\end{equation}

{\bf Remark 1.} (i) The assumption of the theorem means that the
pair $(X,\frac{1}{n}\Sigma)$ has a non-canonical singularity with
the centre at the point $o$. Explicitly, for some exceptional
divisor $R$ over $X$, the centre of which is the point $o$, the
Noether-Fano inequality
$$
\mathop{\rm ord}\nolimits_R\Sigma>n\cdot a(R,X)
$$
holds, where $a(R,X)$ is the discrepancy of $R$ with respect to
$X$.\vspace{0.1cm}

(ii) The self-intersection $Z=(D_1\circ D_2)$ is the
scheme-theoretic intersection of any two general divisors in
$\Sigma$ which is well defined as $\Sigma$ is free from fixed
components.\vspace{0.1cm}

(iii) When $\mathop{\rm mult}_oX=1$, we get the standard
$4n^2$-inequality, see \cite[Chapter 2]{Pukh13z}. For that reason,
we call the inequality (\ref{24.06.2016.1}) the $4n^2$-inequality
as well. The standard $4n^2$-inequality (for the non-singular
case) was first shown in \cite{Pukh00c} on the basis of the
technique developed in \cite{IM}. Later a different proof was
found by Corti \cite{Co00} and various generalizations of the
$4n^2$-inequality were investigated \cite{Ch05c,Pukh10e}, see
\cite[Chapter 2]{Pukh13z} for more details.

Note that in the smooth case (when $\mathop{\rm mult}\nolimits_o
X=1$) the $4n^2$-inequality holds for $\mathop{\rm dim} X\geqslant
3$ without any additional assumptions. This is because the
exceptional divisor of the blow up of the point $o$ on $X$ is just
the projective space, and in the projective space it is very easy
to bound multiplicities in terms of degrees. Unfortunately, it is
not so easy to do so (in the way we need) for hypersurfaces and
complete intersections, which generate the need for additional
assumptions.\vspace{0.1cm}

The author thanks the referees for a number of useful suggestions,
especially for spotting the insufficient lower bound for
$\mathop{\rm dim} X$ in the first version of the
paper.\vspace{0.3cm}


{\bf 2. Generic complete intersection singularities.} The germ
$(X,o)$ is given by a system of $l$ analytic equations
$$
\begin{array}{ccccccc}
0 & = & q_{1,\mu_1} & + & q_{1,\mu_1+1} & + & \dots\\
 & \dots & & & & &\\
0 & = & q_{l,\mu_l} & + & q_{l,\mu_l+1} & + & \dots\\
\end{array}
$$
in ${\mathbb C}^{M+l}$, where
$2\leqslant\mu_1\leqslant\dots\leqslant\mu_l$, $l\geqslant 1$ and
the polynomials $q_{j,i}$ are homogeneous of degree $i$ in the
coordinates $z_1,\dots,z_{M+l}$; the point $o=(0,\dots,0)$ is the
origin. We denote by
$$
\underline{\mu}=(\mu_1,\dots,\mu_l)
$$
the type of the singularity $o\in X$ and set
$$
\mu=\mu_1\cdots\mu_l=\mathop{\rm mult}\nolimits_o X
$$
to be the multiplicity of the point $o$ (assuming the conditions
of general position for the first polynomials
$q_{1,\mu_1},q_{2,\mu_2},\dots,q_{l,\mu_l}$, stated below). Set
also
$$
|\underline{\mu}|=\mu_1+\dots +\mu_l.
$$

Recall that by assumption $M\geqslant l+|\underline{\mu}|+3$. Let
$P\ni o$ be a linear subspace in ${\mathbb C}^{M+l}$ of dimension
$2l+|\underline{\mu}|+3$. Denote by $X_P$ the intersection $X\cap
P$.\vspace{0.1cm}

{\bf Definition 1.} We say that the complete intersection
singularity $(X,o)$ is {\it generic}, if for a general subspace
$P$ of dimension $2l+|\underline{\mu}|+3$ the singularity $o\in
X_P$ is an isolated singularity, $\mathop{\rm
dim}X_P=l+|\underline{\mu}|+3$ and for the blow up
$$
\varphi_P\colon X^+_P\to X_P
$$
of the point $o$, the variety $X^+_P$ is non-singular in
neighborhood of the exceptional divisor $Q_P=\varphi^{-1}_P(o)$,
which is a non-singular complete intersection
$$
Q_P=\{q_{1,\mu_1}=q_{2,\mu_2}=\dots=q_{l,\mu_l}=0\}\subset{\mathbb
P}^{2l+|\underline{\mu}|+2}
$$
of codimension $l$ and type
$\underline{\mu}=(\mu_1,\dots,\mu_l)$.\vspace{0.1cm}

From now on, we assume that the singularity $o\in X$ is generic.
In particular, by Grothendieck's theorem on factoriality
\cite{CL}, $X$ is a factorial variety near the point
$o$.\vspace{0.3cm}


{\bf 3. Start of the proof.} The idea of the proof is as follows.
We use as a model the proof of the standard $4n^2$-inequality by
means of the technique of counting multiplicities as it is given
in \cite[Chapter 2, Section 2.2]{Pukh13z}. First, we observe that
by inversion of adjunction, the existence of a non-canonical
singularity $R$ implies the existence of another singularity $E$
of the same pair $(X,\frac{1}{n}\Sigma)$ which satisfies a {\it
Noether-Fano type} inequality. The latter is somewhat weaker (but
sufficient for our purposes). However, the new singularity $E$ has
the crucial advantage that its centre on the blow up $X^+$ of the
point $o$ has a high dimension. This is done in the present
section.\vspace{0.1cm}

After that, in Section 4 we resolve the singularity $E$ and use
the assumptions on the singular point $o\in X$ to relate the
multiplicities of the system $\Sigma$ and its self-intersection at
the point $o$ with the multiplicities of the strict transforms of
$\Sigma$ and the self-intersection at the ``higher storeys'' of
the resolution, at the centres of the singularity $E$ on those
``higher storeys''.\vspace{0.1cm}

This done, we apply the technique of counting multiplicities in
word for word the same way as in \cite[Chapter 2, Section
2.2]{Pukh13z} and complete the proof.\vspace{0.1cm}

Let us realize this programme.\vspace{0.1cm}

For a general $(2l+|\underline{\mu}|+3)$-subspace $P$ set
$\Sigma_P=\Sigma|_P$ to be the restriction of $\Sigma$ onto $P$.
By inversion of adjunction \cite{Sh93,Kol93}, the pair
$(X_P,\frac{1}{n}\Sigma_P)$ is not canonical (for
$M>l+|\underline{\mu}|+3$, even non-log canonical, but we do not
need that.) Obviously,
$$
Z_P=Z|_P=(Z\circ X_P)
$$
is the self-intersection of the system $\Sigma_P$ and $\mathop{\rm
mult}_oZ=\mathop{\rm mult}_oZ_P$. Therefore, we may (and will)
assume from the beginning that $M=l+|\underline{\mu}|+3$ and so
$P={\mathbb C}^{M+l}$, so that already the original singularity
$o\in X$ is isolated. Now we omit the index $P$ and write
$$
\varphi\,\,\colon X^+\to X
$$
for the blow up of the point $o$ and $Q=\varphi^{-1}(o)$ for the
exceptional divisor, which is a non-singular complete intersection
of type $\underline{\mu}$ in ${\mathbb
P}^{2l+|\underline{\mu}|+2}$.\vspace{0.1cm}

Now let $\Pi\ni o$ be a general linear subspace of dimension
$|\underline{\mu}|+3$. By the symbol $X_{\Pi}$ we denote the
intersection $X\cap\Pi$. Clearly, $o\in X_{\Pi}\subset\Pi={\mathbb
C}^{|\underline{\mu}|+3}$ is an isolated complete intersection
singularity of codimension $l$. Let $\varphi_{\Pi}\,\,\colon
X^+_{\Pi}\to X_{\Pi}$ be the blow up of the point $o$ and
$Q_{\Pi}=\varphi^{-1}_{\Pi}(o)$ the exceptional divisor. Clearly
$Q_{\Pi}\subset{\mathbb P}^{|\underline{\mu}|+2}$ is a
non-singular complete intersection of type $\underline{\mu}$ (and
codimension $l$).\vspace{0.1cm}

Note that by the adjunction formula for the discrepancy we have
the equality $a(Q_{\Pi},X_{\Pi})=2$.\vspace{0.1cm}

For a general divisor $D\in\Sigma$ and its strict transform
$D^+\in\Sigma^+$ on $X^+$ we have
$$
D^+\sim -\nu\, Q
$$
for some positive integer $\nu$ (recall that we consider a local
situation: $o\in X$ is a germ). Obviously, if $\nu>2n$, then
$$
\mathop{\rm mult}\nolimits_oZ\geqslant\nu^2\mu>4n^2\mu
$$
and the $4n^2$-inequality holds. For that reason, from now on we
assume that
$$
\nu\leqslant 2n.
$$
Setting $D_{\Pi}=D|_{X_\Pi}$, we get $D^+_{\Pi}\sim -\nu\,
Q_{\Pi}$. By the inversion of adjunction the pair
$\left(X_{\Pi},\frac{1}{n}D_{\Pi}\right)$ is not log canonical at
the point $o$, the more so not canonical, so for some exceptional
divisor $E_{\Pi}$ over $X_{\Pi}$ the Noether-Fano inequality
$$
\mathop{\rm ord}\nolimits_{E_{\Pi}}\Sigma_{\Pi}>n
a(E_{\Pi},X_{\Pi})
$$
is satisfied. As $\nu\leqslant 2n$ and $a(Q_{\Pi},X_{\Pi})=2$, we
see that $E_{\Pi}\neq Q_{\Pi}$ and $E_{\Pi}$ is a non log
canonical (and so not canonical) singularity of the pair
$$
\left(X^+_{\Pi},\frac{1}{n}D^+_{\Pi}+\frac{(\nu-2n)}{n}Q_{\Pi}\right)
$$
(the more so, of the pair
$\left(X^+_{\Pi},\frac{1}{n}D^+_{\Pi}\right))$. Denote by
$\Delta_{\Pi}\subset Q_{\Pi}$ the centre of $E_{\Pi}$ on
$X^+_{\Pi}$, an irreducible subvariety in $Q_{\Pi}$.\vspace{0.1cm}

{\bf Proposition 1.} {\it If $\mathop{\rm
codim}\,(\Delta_{\Pi}\subset Q_{\Pi})=1$, then the estimate
$$
\mathop{\rm mult}\nolimits_o Z\geqslant 8n^2\mu
$$
holds.}\vspace{0.1cm}

{\bf Proof.} We note that $\mathop{\rm mult}_oZ=\mathop{\rm
mult}_oZ_{\Pi}$. Arguing as in the proof of Proposition 4.1 in
\cite[Chapter 2]{Pukh13z} (see also \cite[Section 1.7]{Ch05c}), we
get the following estimate:
$$
\mathop{\rm
mult}\nolimits_oZ_{\Pi}\geqslant\nu^2\mu+
4\left(3-\frac{\nu}{n}\right)n^2\mu,
$$
and easy calculations complete the proof. Q.E.D.\vspace{0.1cm}

Therefore we may assume that $\mathop{\rm
codim}(\Delta_{\Pi}\subset Q_{\Pi})\geqslant 2$.\vspace{0.1cm}

Coming back to the variety $X$, we conclude that for some
exceptional divisor $E$ over $X$ with the centre at $o$ the
Noether-Fano type inequality
$$
\mathop{\rm ord}\nolimits_E\Sigma>n(2\mathop{\rm
ord}\nolimits_EQ+a(E,X^+))
$$
is satisfied. Moreover, the centre $\Delta\subset Q$ of $E$ on $X$
has codimension at least 2 and dimension at least
$2l$.\vspace{0.3cm}


{\bf 4. Resolution of the singularity $E$.} Consider the sequence
of blow ups
$$
X_0=X\leftarrow X_1=X^+\leftarrow X_2\leftarrow\dots\leftarrow
X_K,
$$
where $\varphi_{i,i-1}\colon X_i\to X_{i-1}$ is the blow up of the
centre $B_{i-1}\subset X_{i-1}$ of the exceptional divisor $E$ on
$X_{i-1}$. In particular, $B_0=o$ and $B_1=\Delta$. Using the
notations, identical to those in \cite[Chapter 2, Section
2.2]{Pukh13z}, we set
$$
E_i=\varphi^{-1}_{i,i-1}(B_{i-1})\subset X_i
$$
to be the exceptional divisor, so that $E_1=Q$. As $X_1=X^+$ is
non-singular in a neighborhood of $E_1$, all subsequent varieties
$X_i$ are non-singular at the generic point of $B_i$ and all
constructions of \cite[Chapter 2, Section 2.2]{Pukh13z} work
automatically for the blow ups $\varphi_{i,i-1}$ with $i\geqslant
2$.\vspace{0.1cm}

The last exceptional divisor $E_K$ defines the discrete valuation
$\mathop{\rm ord}_E$.\vspace{0.1cm}

We divide the sequence $\varphi_{i,i-1}$, $i=1,\dots,K$, of blow
ups into the {\it lower part}, $i=1,\dots,L\leqslant K$,
corresponding to the centres $B_{i-1}$ of codimensions at least 3,
and the {\it upper part}, $i=L+1,\dots,K$, corresponding to the
centres $B_{i-1}$ of codimension two. As usual, we denote the
strict transform of any geometric object on $X_i$ by adding the
upper index $i$ and set:
$$
\nu_i=\mathop{\rm mult}\nolimits_{B_{i-1}}\Sigma^i
$$
for any $i=2,\dots,K$ to be the elementary multiplicities. Let
$\Gamma$ be the oriented graph of the resolution of the
singularity $E$ and $p_{ij}$ the number of paths from the vertex
$i$ to the vertex $j$, $p_{ii}=1$ by definition (see \cite[Chapter
2, Section 2.2]{Pukh13z} for the standard details). We also set
$p_i=p_{Ki}$, $i=1,\dots,K$. Now the Noether-Fano type inequality
takes the form
\begin{equation}\label{27.06.2016.1}
\sum^K_{i=1}p_i\nu_i>\left(2p_1+\sum^K_{i=2}p_i\delta_i\right),
\end{equation}
where $\nu_1=\nu$ and $\delta_i=\mathop{\rm codim}(B_{i-1}\subset
X_{i-1})$ are the elementary discrepancies. By the linearity of
the Noether-Fano type inequality (\ref{27.06.2016.1}) and the
standard properties of the numbers $p_{ij}$ we may assume that
$\nu_K>n$ (replacing, if necessary, $E_K$ by a lower singularity
$E_j$ for some $j<K$). In order to proceed, we need the following
known fact.\vspace{0.1cm}

{\bf Proposition 2.} {\it Let $Y\subset{\mathbb P}^N$ be a
non-singular complete intersection of codimen\-sion $l\geqslant
1$, $S\subset Y$ an irreducible subvariety of codimension
$a\geqslant 1$ and $B\subset Y$ an irreducible subvariety of
dimension $al$, where the estimate $N\geqslant(l+1)(a+1)$ is
satisfied. Then the inequality
$$
\mathop{\rm mult}\nolimits_BS\leqslant m
$$
holds, where $m\geqslant 1$ is defined by the condition $S\sim
mH^a_Y$ and $H_Y\in A^1Y$ is the class of a hyperplane section of
$Y$.}\vspace{0.1cm}

{\bf Proof} for the case $l=1$ was given in \cite{Pukh02b}. The
argument extends directly to the general case of an arbitrary $l$,
see \cite{Suzuki15} (also \cite{Pukh06b,Ch05a}).
Q.E.D.\vspace{0.1cm}

Applying Proposition 2 to a divisor in the linear system
$\Sigma^1|_Q$, we conclude that $\nu_1\geqslant\nu_2$, since
$\mathop{\rm dim}B_1=\mathop{\rm dim}\Delta\geqslant 2l$. The
inequalities
$$
\nu_2\geqslant\nu_3\geqslant\dots\geqslant\nu_K
$$
are standard. We deduce that the upper part of the resolution of
$E$ is non-empty (that is to say, $L<K$) and the upper part of the
graph $\Gamma$ is a chain:
$$
L\leftarrow(L+1)\leftarrow\dots\leftarrow K;
$$
moreover, there are no arrows connecting either of the vertices
$L+1,\dots,K$ with any of vertices $1,2,\dots,L-1$. (These are the
standard consequences of inequalities $\nu_K>n$ and
$\nu_1\leqslant 2n$, see \cite[Chapter 2, Section 2.2]{Pukh13z}.)
We do not need this additional information for the proof of our
theorem, but in particular geometric problems it might be
useful.\vspace{0.3cm}


{\bf 5. The technique of counting multiplicities.} Now everything
is ready for the proof of the desired inequality
(\ref{24.06.2016.1}). Take a general pair of divisors
$D_1,D_2\in\Sigma$ and set
$$
Z=Z_0=(D_1\circ D_2)
$$
to be their scheme-theoretic intersection, the self-intersection
of the mobile linear system $\Sigma$. Recall that the strict
transform of an irreducible subvariety or an effective cycle, or a
linear system on some $X_i$ is denoted by adding the upper index
$i$. (This notation silently implies that the irreducible
subvariety or the effective cycle etc. is sitting on a lower
storey $X_j$, $j\leqslant i$, of the resolution and that the
operation of taking the strict transform is well defined for that
particular subvariety etc.) For $i\geqslant 1$ write
$$
(D^i_1\circ D^i_2)=(D^{i-1}_1\circ D^{i-1}_2)^i+Z_i,
$$
where the effective cycle $Z_i$ of codimension 2 is supported on
$E_i$ and so may be viewed as an effective divisor on $E_i$. Thus
for any $i\leqslant L$ we obtain the presentation
$$
(D^i_1\circ D^i_2)=Z^i_0+Z^i_1+\dots+Z^i_{i-1}+Z_i.
$$
For any $j>i$, where $j\leqslant L$, set
$$
m_{i,j}=\mathop{\rm mult}\nolimits_{B_{j-1}}Z^{j-1}_i
$$
and for $i=2,\dots,L$ set $d_i=\mathop{\rm deg}Z_i$ in the same
sense as in \cite[Chapter 2, Section 2.2]{Pukh13z}. For the
effective divisor $Z_1$ on $E_1=Q$ we have the relation
$$
Z_1\sim d_1H_Q
$$
for some $d_1\in{\mathbb Z}_+$, where $H_Q$ is the class of a
hyperplane section of the complete intersection $Q\subset{\mathbb
P}^{4l+2}$. Following the procedure of \cite[Chapter 2]{Pukh13z},
we obtain the system of equalities
$$
\begin{array}{rccccccc}
\mu\!\!\!\!\! & (\nu_1^2+d_1) & = & m_{0,1}, & & & &  \\
    &  \nu_2^2+d_2  & = & m_{0,2} & + & m_{1,2}, &  & \\
    &               & \dots &  &  &  &  & \\
    &  \nu_i^2+d_i  & = & m_{0,i} & + & \dots & + & m_{i-1,i}, \\
    &               & \dots &  &  &  &  &
\end{array}
$$
$i=2,\dots,L$, where the estimate
$$
d_L\geqslant\sum^K_{i=L+1}\nu^2_i
$$
holds as usual, see \cite[p. 53]{Pukh13z}.\vspace{0.1cm}

{\bf Proposition 3.} (i) {\it The inequality
$$
d_1\geqslant m_{1,2}
$$
holds.}\vspace{0.1cm}

(ii) {\it The inequality
$$
m_{0,1}\geqslant\mu m_{0,2}
$$
holds.}\vspace{0.1cm}

{\bf Proof.} Part (i) follows from Proposition 2 as $Z_1\sim d_1
H_Q$ and $\mathop{\rm dim}B_1\geqslant 2l$. In order to show part
(ii), we note that (numerically)
$$
(Z^1\circ E_1)\sim\frac{1}{\mu}m_{0,1}H^2_Q
$$
as $m_{0,1}=\mathop{\rm deg}(Z^1\circ E_1)$, the cycle $(Z^1\circ
E_1)=(Z^1\circ Q)$ being of pure codimension 2 on $Q$. Applying
Proposition 2 to the cycle $(Z^1\circ Q)$, we get the inequality
$$
m_{0,2}\leqslant\mathop{\rm mult}\nolimits_{\Delta}(Z^1\circ
Q)\leqslant\frac{1}{\mu}m_{0,1},
$$
which completes the proof of the proposition. Q.E.D.\vspace{0.1cm}

The more so, $m_{0,1}\geqslant\mu m_{0,i}$ for $i\geqslant 3$ as
$m_{0,2}\geqslant m_{0,3}\geqslant\dots\geqslant
m_{0,L}$.\vspace{0.1cm}

Now set
$$
m^*_{i,j}=\mu m_{i,j}
$$
for $(i,j)\neq (0,1)$ and $m^*_{0,1}=m_{0,1}$. Also set
$$
d^*_i=\mu d_i
$$
for $i=1,\dots,L$. We obtain the following system of equalities:
$$
\begin{array}{ccccccc}
\mu\nu_1^2+d^*_1 & = & m^*_{0,1}, & & & &  \\
\mu\nu_2^2+d^*_2  & = & m^*_{0,2} & + & m^*_{1,2}, &  & \\
               & \dots &  &  &  &  & \\
\mu\nu_i^2+d^*_i  & = & m^*_{0,i} & + & \dots & + & m^*_{i-1,i}, \\
               & \dots &  &  &  &  &
\end{array}
$$
where $i=1,\dots,L$, and
$$
d^*_L\geqslant\mu\sum^K_{i=L+1}\nu^2_i,
$$
where the integers $m^*_{i,j}$ and $d^*_i$ satisfy precisely the
same properties, as the integers $m_{i,j}$ and $d_i$ in the
non-singular case considered in \cite[Chapter 2, p.
52-53]{Pukh13z}. Now repeating the arguments of \cite[Chapter 2,
p. 52-53]{Pukh13z} word for word, we obtain the inequality
$$
\left(\sum^L_{i=1}p_i\right)\mathop{\rm
mult}\nolimits_oZ\geqslant\mu\sum^K_{i=1}p_i\nu^2_i,
$$
which in the standard way implies the desired estimate
$$
\mathop{\rm mult}\nolimits_oZ>\mu\cdot 4n^2.
$$
Proof of the theorem is completed.\vspace{0.1cm}

{\bf Remark 2.} The inequality (\ref{24.06.2016.1}) essentially
simplifies the proof of birational superrigidity of Fano
hypersurfaces with isolated singularities of general position
given in \cite{Pukh02a}. The cases of singular points of
multiplicity $\mu=3$ and 4 in that paper are really hard. The
inequality (\ref{24.06.2016.1}) gives for the multiplicity
$\mathop{\rm mult}_oZ$ at such points the lower bound $12n^2$ and
$16n^2$, respectively, which is more than enough to exclude the
maximal singularities over such points by the standard (in fact,
relaxed) technique of hypertangent divisors. More applications of
the inequality (\ref{24.06.2016.1}) in the spirit of \cite{EP,EvP}
will be given separately.

\begin{flushleft}
Department of Mathematical Sciences,\\
The University of Liverpool
\end{flushleft}

\noindent{\it pukh@liverpool.ac.uk}

\end{document}